\documentclass{elsart}

\usepackage{amssymb}
\usepackage{amsmath}
\usepackage{graphicx}
\usepackage{algorithm}
\usepackage{algorithmic}
\usepackage{url}

\newcommand{\PN}{2^N}
\newcommand{\R}{\mathbb{R}}                     
\newcommand{\Cdot}{\boldsymbol{\cdot}}

\theoremstyle{plain}
\newtheorem{exm}{Example}

\begin{document}

\begin{frontmatter}

\title{On the moments and distribution of discrete Choquet integrals from continuous distributions}

\author[ivan]{Ivan Kojadinovic},
\ead{ivan@stat.auckland.ac.nz}
\author[jean-luc]{Jean-Luc Marichal}
\ead{jean-luc.marichal@uni.lu}
\address[ivan]{Department of Statistics, The University of Auckland, Private Bag 92019, Auckland 1142, New Zealand.}
\address[jean-luc]{Mathematics Research Unit, University of Luxembourg, 162A, avenue de la Fa\"iencerie, L-1511 Luxembourg, G.D. Luxembourg.}

\date{Revised version}

\begin{abstract}
We study the moments and the distribution of the discrete Choquet integral when regarded as a real function of a random sample drawn from a
continuous distribution. Since the discrete Choquet integral includes weighted arithmetic means, ordered weighted averaging functions, and
lattice polynomial functions as particular cases, our results encompass the corresponding results for these aggregation functions. After
detailing the results obtained in \cite{MarKoj08} in the uniform case, we present results for the standard exponential case, show how
approximations of the moments can be obtained for other continuous distributions such as the standard normal, and elaborate on the asymptotic
distribution of the Choquet integral. The results presented in this work can be used to improve the interpretation of discrete Choquet integrals
when employed as aggregation functions.
\end{abstract}


\end{frontmatter}

\section{Introduction}

Aggregation functions are of central importance in many fields such as statistics, information fusion, risk analysis, or decision theory. In
this paper, the primary object of interest is a natural extension of the weighted arithmetic mean known as the (discrete) {\em Choquet integral}
\cite{Cho53,Gra95,Mar00}. Also known in discrete mathematics as the {\em Lov\'asz extension} of pseudo-Boolean functions \cite{Lov83}, the
Choquet integral is a very flexible aggregation function that includes weighted arithmetic means, ordered weighted averaging
functions~\cite{Yag88}, and lattice polynomial functions as special cases~\cite{Mar06,MarKoj08}.

Although the Choquet integral has been extensively employed as an aggregation function (see e.g.\ \cite{GraMurSug00} for an overview), its
moments and its distribution seem to have never been thoroughly studied from a theoretical perspective. The aim of this work is to attempt to
fill this gap in the case when the Choquet integral is regarded as a real function of a random sample drawn from a continuous distribution.

The starting point of our study is a natural distributional relationship between linear combinations of order statistics and the Choquet
integral, which merely results from the piecewise linear decomposition of the latter. As a consequence, exact formulations of the moments and
the distribution of the Choquet integral can be provided whenever exact formulations are known for linear combinations of order statistics.
Likewise, approximation and asymptotic results can be provided whenever available for linear combinations of order statistics.

The paper is organized as follows. In the second section, we recall the definition of the discrete Choquet integral. The third section is
devoted to the expression of the distribution (resp.\ the moments) of the Choquet integral in terms of the distribution (resp.\ the moments) of
linear combinations of order statistics. The case of standard uniform input variables is treated in Section~\ref{sec:uniform}. More precisely,
the results obtained in~\cite{MarKoj08} are detailed, and algorithms for computing the probability density function (p.d.f.) and the cumulative
distribution function (c.d.f.) of the Choquet integral are provided. The fifth section deals with the standard exponential case, while the sixth
one shows how approximations of moments can be obtained for other continuous distribution such as the standard normal. In the last section, we
discuss conditions under which the asymptotic distribution of the Choquet integral is a mixture of normals.

The results obtained in this work have numerous applications. The most immediate ones are related to the interpretation of the Choquet integral
when seen as an aggregation function. In multicriteria decision aiding in particular, the presented results can be used to generalize the
behavioral indices studied e.g. in \cite{Mar00b,Mar04}. In classifier fusion, they can enable a theoretical study of the so-called {\em fuzzy
approach} to classifier combination (see e.g.\ \cite{Kun03}) in the spirit of that done in~\cite{Kun02}.

Note that most of the methods and algorithms discussed in this work have been implemented in the R package {\tt kappalab} \cite{kappalab}
available on the Comprehensive R Archive Network (\url{http://CRAN.R-project.org}).

\section{The discrete Choquet integral}

Define $N:=\{1,\dots,n\}$ as a set of attributes, criteria, or players, and denote by $\frak{S}_n$ the set of permutations on $N$. A set
function $\nu:\PN \to [0,1]$ is said to be a {\em game} on $N$ if it satisfies $\nu(\varnothing) = 0$.

\begin{defn}
\label{def_Choquet} The discrete Choquet integral of $\textbf{x} \in \mathbb{R}^n$ w.r.t.\ a game $\nu$ on $N$ is defined by
$$
C_\nu(\textbf{x}) := \sum_{i=1}^n p^{\nu,\sigma}_i \, x_{\sigma(i)} \, ,
$$
where $\sigma\in\frak{S}_n$ is such that $x_{\sigma(1)} \geqslant \dots \geqslant x_{\sigma(n)}$, where
$$
p^{\nu,\sigma}_i := \nu^\sigma_i - \nu^\sigma_{i-1}, \qquad \forall \, i \in N,
$$
and where $\nu^\sigma_i := \nu\big(\{\sigma(1),\dots,\sigma(i)\}\big)$ for any $i=0,\dots,n$. In particular, $\nu^\sigma_{0} := 0$.
\end{defn}

Note that the permutation $\sigma$ in the defintion of the Choquet integral of $\textbf{x}$ is traditionally taken such that $x_{\sigma(1)}
\leqslant \dots \leqslant x_{\sigma(n)}$. The reason for not adopting this convention in this work is due to the fact that it would have led to
much more complicated expressions of the results to be presented in Section~\ref{sec:uniform}.

From the above definition, we see that the Choquet integral is a piecewise linear function that coincides with a weighted sum on each
$n$-dimensional polyhedron
\begin{equation}\label{eq:Rsigma}
R_{\sigma}:=\{\textbf{x}\in \R^n\mid x_{\sigma(1)}\geqslant \cdots \geqslant x_{\sigma(n)}\},\qquad \sigma\in\frak{S}_n,
\end{equation}
whose union covers $\R^n$. It can additionally be immediately verified that it is a continuous function.

When defined as above, the Choquet integral coincides with the {\em Lov\'asz extension} \cite{Lov83} of the unique pseudo-Boolean function that
can be associated with $\nu$ \cite{Mar00d} and can be alternatively regarded as a linear combination of lattice polynomial functions (see e.g.\
\cite{MarKoj08}).

In aggregation theory, it is natural to additionally require that the game $\nu$ is monotone w.r.t.\ inclusion and satisfies $\nu(N) = 1$, in
which case it is called a {\em capacity} \cite{Cho53}. The resulting aggregation function $C_\nu$ is then nondecreasing in each variable and
coincides with a weighted arithmetic mean on each of the $n$-dimensional polyhedra defined by~(\ref{eq:Rsigma}). Furthermore, in this case, for
any $T \subseteq N$, the coefficient $\nu(T)$ can be naturally interpreted as the {\em weight} or the {\em importance} of the subset $T$ of
attributes \cite{Mar00}.

The Choquet integral w.r.t.\ a capacity satisfies very appealing properties for aggregation. For instance, it is comprised between the minimum
and the maximum, stable under the same transformations of interval scales in the sense of the theory of measurement, and coincides with a
weighted arithmetic mean whenever the capacity is additive. An axiomatic characterization is provided in~\cite{Mar00}. Moreover, the Choquet
integral w.r.t.\ a capacity includes weighted arithmetic means, ordered weighted averaging functions \cite{Yag88}, and lattice polynomial
functions as particular cases \cite{Mar06,MarKoj08}.

\section{Distributional relationships with linear combinations of order statistics}

In the present section, we investigate the moments and the distribution of the Choquet integral when considered as a function of $n$ continuous
i.i.d.\ random variables. Our main theoretical results, stated in the following proposition and its corollary, yield expressions of the moments
and the distribution of the Choquet integral in terms of the moments and the distribution of linear combinations of order statistics.

Let $X_1,\dots,X_n$ be a random sample drawn from a continuous c.d.f.\ $F:\R\to\R$ with associated p.d.f.\ $f:\R\to\R$, and let $X_{1:n}
\leqslant \cdots \leqslant X_{n:n}$ denote the corresponding order statistics. Furthermore, let
\begin{eqnarray*}
Y_\nu &:=& C_\nu(X_1,\dots,X_n),\\
Y_\nu^\sigma &:=& \sum_{i=1}^n p^{\nu,\sigma}_i \, X_{n-i+1:n},\qquad \sigma \in \frak{S}_n.
\end{eqnarray*}
Let also $F_{\nu}(y)$ and $F_{\nu}^\sigma(y)$ be the c.d.f.s of $Y_\nu$ and $Y_\nu^\sigma$, respectively. Finally, let $h:\R\to\R$ be any
measurable function.

\begin{prop}\label{prop:main}
For any game $\nu$ on $N$, we have
$$
\mathbf{E}[h(Y_\nu)] = \frac{1}{n!} \sum_{\sigma \in \frak{S}_n} \mathbf{E} [h(Y_\nu^\sigma)].
$$
\end{prop}

\begin{pf*}{Proof.}
By definition, we have
\begin{align*}
\mathbf{E}[h(Y_\nu)] &= \int_{\R^n} h(C_\nu(x_1,\dots,x_n)) \prod_{i=1}^n f(x_i)\,\mathrm{d} x_i \\
&= \sum_{\sigma \in \frak{S}_n} \int_{R_\sigma} h \Big( \sum_{i=1}^n p^{\nu,\sigma}_i \, x_{\sigma(i)} \Big) \prod_{i=1}^n f(x_i)
\,\mathrm{d}x_i .
\end{align*}
Using the well-known fact (see e.g.\ \cite[\S2.2]{DavNag03}) that the joint p.d.f.\ of $X_{1:n} \leqslant \cdots \leqslant X_{n:n}$ is
$$
n! \prod_{i=1}^n f(x_i), \qquad x_1 \leqslant \cdots \leqslant x_n,
$$
we obtain
$$
\mathbf{E}[h(Y_\nu)] = \frac{1}{n!} \sum_{\sigma \in \frak{S}_n} \mathbf{E} \Big[ h \Big( \sum_{i=1}^n p^{\nu,\sigma}_i \, X_{n-i+1:n} \Big)
\Big],
$$
which completes the proof. \qed \end{pf*}

Before going through the main corollary, recall that the {\em plus}\/ (resp.\ {\em minus}) {\em truncated power function}\/ $x^n_+$ (resp.\
$x^n_-$) is defined to be $x^n$ if $x>0$ (resp.\ $x<0$) and zero otherwise.

\begin{cor}\label{cor:main}
For any game $\nu$ on $N$, we have
$$
F_\nu(y) = \frac{1}{n!} \sum_{\sigma \in \frak{S}_n} F_\nu^\sigma(y).
$$
\end{cor}

\begin{pf*}{Proof.}
Define $h_y(x):=(x-y)_-^0$. Then, from Proposition~\ref{prop:main}, for any $y\in\R$, we have
$$
F_\nu(y)=\mathbf{E}[h_y(Y_\nu)]=\frac{1}{n!} \sum_{\sigma \in \frak{S}_n} \mathbf{E} [h_y(Y_\nu^\sigma)]=\frac{1}{n!} \sum_{\sigma \in
\frak{S}_n} F_\nu^\sigma(y).
$$
\qed \end{pf*}

The results stated in Proposition~\ref{prop:main} and Corollary~\ref{cor:main} are not very surprising. From Definition~\ref{def_Choquet}, it is
clear that the Choquet integral is a linear combination of order statistics whose coefficients depend on the ordering of the arguments. The
different possible orderings merely lead to a division of the integration domain $\R^n$ into the subdomains $R_\sigma$ ($\sigma \in \frak{S}_n$)
defined in~(\ref{eq:Rsigma}), and the difficult part still lies in the evaluation of the moments and the distribution of linear combinations of
order statistics.

The relationship for the raw moments is obtained by considering the special case $h(x)=x^r$, which may still lead to tedious computations. From
Proposition~\ref{prop:main}, we obtain
$$
\mathbf{E}[Y_\nu] = \frac{1}{n!} \sum_{\sigma \in \frak{S}_n}\sum_{k=1}^n p_k^{\nu,\sigma}\,\mathbf{E}[X_{n-k+1:n}],
$$
and more generally,
$$
\mathbf{E}[Y_\nu^r] = \frac{1}{n!} \sum_{\sigma \in \frak{S}_n}\sum_{k_1,\ldots,k_r=1}^n\Big(\prod_{i=1}^r
p_{k_i}^{\nu,\sigma}\Big)\,\mathbf{E}\Big[\prod_{i=1}^r X_{n-k_i+1:n}\Big].
$$
Unfortunately, this latter formula involves a huge number of terms, namely $n!\, n^r$.
The following result (see \cite[Prop.~3]{MarKoj08} for the uniform case) yields the $r$th raw moment as a sum of $(r+1)^n$ terms, each of which
is a product of coefficients $\nu(T)$.

\begin{prop}\label{prop:RawMo}
For any integer $r\geqslant 1$ and any game $\nu$ on $N$, setting $T_{r+1}:=N$ and $X_{0:n}:=0$, we have
$$
\mathbf{E}[Y_\nu^r] = \sum_{T_1 \subseteq \cdots \subseteq T_r \subseteq N}\frac{r!}{[T]_0!\cdots
[T]_n!}\,\Big(\prod_{i=1}^r\frac{\nu(T_i)}{{|T_{i+1}|\choose |T_i|}}\Big)\,\mathbf{E}\Big[\prod_{i=1}^r(X_{n-|T_i|+1:n}-X_{n-|T_i|:n})\Big],
$$
where $[T]_j$ represents the number of ``$j$'' among $|T_1|,\ldots,|T_r|$.
\end{prop}

\begin{pf*}{Proof.}
Fix $\sigma \in \frak{S}_n$. Rewriting $Y_\nu^\sigma$ as
$$
Y_\nu^\sigma = \sum_{i=0}^n \nu_i^{\sigma} (X_{n-i+1:n}-X_{n-i:n}),
$$
and then using the multinomial theorem, we obtain
\begin{eqnarray*}
(Y_\nu^\sigma)^r &=& \sum_{\textstyle{r_1,\ldots,r_n\geqslant 0\atop r_1+\cdots +r_n=r}} \frac{r!}{r_1!\cdots r_n!} \,\prod_{i=0}^n
(\nu_i^{\sigma})^{r_i}\, (X_{n-i+1:n}-X_{n-i:n})^{r_i}\\
&=& \sum_{0\leqslant i_1\leqslant\cdots\leqslant i_r\leqslant n} \frac{r!}{[i]_0!\cdots [i]_n!} \,\prod_{k=1}^r \nu_{i_k}^{\sigma}\,
(X_{n-i_k+1:n}-X_{n-i_k:n}),
\end{eqnarray*}
where $[i]_j$ represents the number of ``$j$'' among $i_1,\ldots,i_r$. Now, using Proposition~\ref{prop:main} with $h(x)=x^r$, we immediately
obtain
$$
\mathbf{E}[Y_\nu^r] = \sum_{0\leqslant i_1\leqslant\cdots\leqslant i_r\leqslant n} \frac{r!}{[i]_0!\cdots [i]_n!}\, \mathbf{E}\Big[\prod_{k=1}^r
(X_{n-i_k+1:n}-X_{n-i_k:n})\Big] \, \frac{1}{n!} \sum_{\sigma \in \frak{S}_n} \prod_{k=1}^r \nu_{i_k}^{\sigma}.
$$
The final result then follows from the identity (see the proof of \cite[Prop.~3]{MarKoj08})
$$
\frac{1}{n!}\sum_{\sigma\in\frak{S}_n}\prod_{k=1}^r \nu_{i_k}^{\sigma} = \sum_{\textstyle{T_1 \subseteq \cdots \subseteq T_r \subseteq N \atop
|T_1|=i_1, \dots, |T_r|=i_r}}\,\prod_{i=1}^r\frac{\nu(T_i)}{{|T_{i+1}|\choose |T_i|}}.
$$
\qed \end{pf*}

For example, the first two raw moments are
\begin{equation}
\label{eq:expect} \mathbf{E}[Y_\nu] = \sum_{T\subseteq N} \frac{\nu(T)}{{n\choose |T|}}\,\mathbf{E}[X_{n-|T|+1:n}-X_{n-|T|:n}]
\end{equation}
and
$$
\mathbf{E}[Y_\nu^2] = \sum_{T_1 \subseteq T_2 \subseteq N}\frac{2}{[T]_0!\cdots [T]_n!}\,\frac{\nu(T_1)\nu(T_2)}{{|T_2|\choose |T_1|}{n\choose
|T_2|}}\,\mathbf{E}\big[(X_{n-|T_1|+1:n}-X_{n-|T_1|:n})(X_{n-|T_2|+1:n}-X_{n-|T_2|:n})\big],
$$
that is,
\begin{multline}
\label{eq:expect_square}
\mathbf{E}[Y_\nu^2] = \sum_{T_1 \varsubsetneq T_2 \subseteq N}2\,\frac{\nu(T_1)\nu(T_2)}{{|T_2|\choose |T_1|}{n\choose |T_2|}}\,\mathbf{E}\big[(X_{n-|T_1|+1:n}-X_{n-|T_1|:n})(X_{n-|T_2|+1:n}-X_{n-|T_2|:n})\big]\\
\null + \sum_{T \subseteq N}\frac{\nu(T)^2}{{n\choose |T|}}\,\mathbf{E}\big[X_{n-|T|+1:n}-X_{n-|T|:n}\big]^2.
\end{multline}

\section{The uniform case}
\label{sec:uniform}

In this section, we focus on the moments and the distribution of $Y_\nu$ when the random sample $X_1,\dots,X_n$ is drawn from the standard
uniform distribution. To emphasize this last point, as classically done, we shall denote the random sample as $U_1,\dots,U_n$ and the
corresponding order statistics by $U_{1:n} \leqslant \dots \leqslant U_{n:n}$.

Before detailing the results obtained in \cite{MarKoj08} and providing algorithms for computing the p.d.f.\ and the c.d.f.\ of the Choquet
integral, we recall some basic material related to divided differences (see e.g.\ \cite{Dav75,DeVLor93,Pow81} for further details).

\subsection{Divided differences}

Let $\mathcal{A}^{(n)}$ be the set of $n-1$ times differentiable one-place functions $g$ such that $g^{(n-1)}$ is absolutely continuous. The
$n$th {\em divided difference}\/ of a function $g\in \mathcal{A}^{(n)}$ is the symmetric function of $n+1$ arguments defined inductively by
$\Delta[g:a_0]:=g(a_0)$ and
$$
\Delta[g:a_0,\ldots,a_n]:=
\begin{cases}
\displaystyle{\frac{\Delta[g:a_1,\ldots,a_n]-\Delta[g:a_0,\ldots,a_{n-1}]}{a_n-a_0}}, & \mbox{if $a_0\neq a_n$},\\
\displaystyle{\frac{\partial}{\partial a_0}\,\Delta[g:a_0,\ldots,a_{n-1}]}, & \mbox{if $a_0= a_n$.}
\end{cases}
$$

The {\em Peano representation}\/ of the divided differences is given by 
$$
\Delta[g:a_0,\ldots,a_n]=\frac 1{n!}\,\int_{\R} g^{(n)}(t)\, M(t\mid a_0,\ldots,a_n)\,\mathrm{d}t,
$$
where $M(t\mid a_0,\ldots,a_n)$ is the {\em B-spline} of order $n$, with knots $\{a_0,\ldots,a_n\}$, defined as
\begin{equation}\label{eq:BSpline}
M(t\mid a_0,\ldots,a_n):=n\,\Delta[(\Cdot-t)^{n-1}_+:a_0,\ldots,a_n].
\end{equation}

We also recall the {\em Hermite-Genocchi formula}: For any function $g\in \mathcal{A}^{(n)}$, we have
\begin{equation}\label{eq:HermGeno}
\Delta[g:a_0,\ldots,a_n]=
\int_{R_{id} \cap [0,1]^n} g^{(n)}\Big[a_0+\sum_{i=1}^n (a_i-a_{i-1})x_i\Big]\,\mathrm{d}x,
\end{equation}
where $R_{id}$ is the region defined in (\ref{eq:Rsigma}) when $\sigma$ is the identity permutation.

For distinct arguments $a_0,\ldots,a_n$, we also have the following formula, which can be verified by induction,
\begin{equation}\label{eq:DvDDistArg}
\Delta[g:a_0,\ldots,a_n]=\sum_{i=0}^n\frac{g(a_i)}{\prod_{j\neq i}(a_i-a_j)}.
\end{equation}

\subsection{Moments and distribution}

Let $g\in \mathcal{A}^{(n)}$. From~(\ref{eq:HermGeno}), we immediately have that
\begin{equation}\label{eq:rztt}
\mathbf{E} \Big[ g^{(n)} \Big(\sum_{i=1}^n p^{\nu,\sigma}_i \, U_{n-i+1:n} \Big) \Big] = n!\, \Delta [g : \nu_0^\sigma, \dots, \nu_n^\sigma]
\end{equation}
since the joint p.d.f.\ of $U_{1:n} \leqslant \dots \leqslant U_{n:n}$ is $1/n!$ on $R_{id} \cap [0,1]^n$ and zero elsewhere.

Now, combining (\ref{eq:rztt}) with Proposition~\ref{prop:main}, we obtain
\begin{equation}
\label{eq:ExpgY} \mathbf{E}[g^{(n)}(Y_\nu)]=\sum_{\sigma\in\frak{S}_n}\Delta[g:\nu_0^{\sigma},\dots,\nu_n^{\sigma}].
\end{equation}
Eq.~(\ref{eq:ExpgY}) provides the expectation $\mathbf{E}[g^{(n)}(Y_\nu)]$ in terms of the divided differences of $g$ with arguments
$\nu_0^{\sigma},\ldots,\nu_n^{\sigma}$ $(\sigma\in\frak{S}_n)$. An explicit formula can be obtained by (\ref{eq:DvDDistArg}) whenever the
arguments are distinct for every $\sigma\in\frak{S}_n$.

Clearly, the special cases
$$
g(x)=\frac{r!}{(n+r)!}\, x^{n+r},~\frac{r!}{(n+r)!}\, [x-\mathbf{E}(Y_\nu)]^{n+r},~\mbox{and}~\frac{e^{tx}}{t^n}
$$
give, respectively, the raw moments, the central moments, and the moment-generating function of $Y_\nu$. As far as the raw moments are
concerned, we have the following result \cite[Prop.~3]{MarKoj08}, which is a special case of Proposition~\ref{prop:RawMo}.

\begin{prop}\label{prop:RM}
For any integer $r\geqslant 1$ and any game $\nu$ on $N$, setting $T_{r+1}:=N$, we have
$$
\mathbf{E}[Y_\nu^r] = \frac{1}{{n+r\choose r}}\,\sum_{T_1 \subseteq \cdots \subseteq T_r \subseteq
N}~\prod_{i=1}^r\frac{\nu(T_i)}{{|T_{i+1}|\choose |T_i|}}.
$$
\end{prop}

Proposition~\ref{prop:RM} provides an explicit expression for the $r$th raw moment of $Y_\nu$ as a sum of $(r+1)^n$ terms. For instance, the
first two moments are
\begin{eqnarray}
\label{ECh_unif} \mathbf{E}[Y_\nu] &=&\frac{1}{n+1}\,\sum_{T\subseteq N}\frac{\nu(T)}{{n\choose |T|}},\\
\label{ECh2_unif} \mathbf{E}[Y_\nu^2] &=& \frac{2}{(n+1)(n+2)}\,\sum_{T_1\subseteq T_2\subseteq N}\frac{\nu(T_1)\nu(T_2)}{{|T_2|\choose
|T_1|}{n\choose |T_2|}}.
\end{eqnarray}

By using (\ref{eq:ExpgY}) with $g(x)=\frac{1}{n!} (x-y)^n_-$, we also obtain the c.d.f.\ $F_\nu(y)$ of $Y_\nu$ \cite{MarKoj08}.

\begin{thm}\label{thm:CDFYn}
There holds
\begin{equation}\label{eq:CDFYnu}
F_\nu(y)= \frac {1}{n!}\,\sum_{\sigma\in\frak{S}_n}\Delta[(\Cdot -y)^n_-:\nu_0^{\sigma},\ldots,\nu_n^{\sigma}] = 1-\frac
{1}{n!}\,\sum_{\sigma\in\frak{S}_n}\Delta[(\Cdot -y)^n_+:\nu_0^{\sigma},\ldots,\nu_n^{\sigma}].
\end{equation}
\end{thm}

It follows from (\ref{eq:CDFYnu}) that the p.d.f.\ of $Y_\nu$ is simply given by
\begin{align}\label{eq:fDens}
\nonumber
f_\nu(y) &= -\frac{1}{(n-1)!}\,\sum_{\sigma\in\frak{S}_n}\Delta[(\Cdot -y)^{n-1}_-:\nu_0^{\sigma},\ldots,\nu_n^{\sigma}] \\
&=\frac {1}{(n-1)!}\,\sum_{\sigma\in\frak{S}_n}\Delta[(\Cdot -y)^{n-1}_+:\nu_0^{\sigma},\ldots,\nu_n^{\sigma}],
\end{align}
or, using the B-spline notation (\ref{eq:BSpline}), by
$$
f_\nu(y)=\frac {1}{n!}\, \sum_{\sigma\in\frak{S}_n} M(y\mid \nu_0^{\sigma},\ldots,\nu_n^{\sigma}).
$$

{\em Remark:}
\begin{enumerate}

\item[(i)] When the arguments $\nu_0^{\sigma},\ldots,\nu_n^{\sigma}$ are distinct for every $\sigma\in\frak{S}_n$, then combining
(\ref{eq:DvDDistArg}) with (\ref{eq:CDFYnu}) immediately yields the following explicit expressions
$$
F_\nu(y)=\frac{1}{n!}\,\sum_{\sigma\in\frak{S}_n}\sum_{i=0}^n\frac{(\nu_i^{\sigma}-y)^n_-}{\prod_{j\neq
i}(\nu_i^{\sigma}-\nu_j^{\sigma})}=1-\frac{1}{n!}\,\sum_{\sigma\in\frak{S}_n}\sum_{i=0}^n\frac{(\nu_i^{\sigma}-y)^n_+}{\prod_{j\neq
i}(\nu_i^{\sigma}-\nu_j^{\sigma})}.
$$


\item[(ii)] The case of linear combinations of order statistics, called {\em ordered weighted averaging operators} in aggregation theory (see
e.g.\ \cite{Yag88}), is of particular interest. In this case, each $\nu_i^{\sigma}$ is independent of $\sigma$, so that we can write
$\nu_i:=\nu_i^{\sigma}$. The main formulas then reduce to (see e.g.\ \cite{AdeSan06,AgaDalSin02})
\begin{eqnarray*}
\mathbf{E}[g^{(n)}(Y_\nu)] &=& n!\, \Delta[g:\nu_0,\ldots,\nu_n],\\
F_\nu(y) &=& \Delta[(\Cdot -y)^n_-:\nu_0,\ldots,\nu_n],\\
f_\nu(y) &=& M(y\mid \nu_0,\ldots,\nu_n).
\end{eqnarray*}
Note also that the Hermite-Genocchi formula (\ref{eq:HermGeno}) provides nice geometric interpretations of $F_\nu(y)$ and $f_\nu(y)$ in terms of
volumes of slices and sections of canonical simplices (see also \cite{Ali73,Ger81}).
\end{enumerate}

\subsection{Algorithms}

Both the functions $F_\nu$ and $f_\nu$ require the computation of divided differences of truncated power functions. On this issue, we recall a
recurrence equation, due to de~Boor \cite{deB72} and rediscovered independently by Varsi \cite{Var73} (see also \cite{Ali73}), which allows to
compute $\Delta[(\Cdot -y)^{n-1}_+:a_0,\ldots,a_n]$ in $O(n^2)$ operations.

Rename as $b_1,\dots,b_r$ the elements $a_i$ such that $a_i < y$ and as $c_1,\dots,c_s$ the elements $a_i$ such that $a_i \geqslant y$ so that
$r+s=n+1$. Then, the unique solution of the recurrence equation
$$
\alpha_{k,l} = \frac{(c_l - y) \alpha_{k-1,l} + (y - b_k) \alpha_{k,l-1}}{c_l - b_k}, \qquad k\leqslant r,\, l\leqslant s,
$$
with initial values $\alpha_{1,1} = (c_1 - b_1)^{-1}$ and $\alpha_{0,l} = \alpha_{k,0} = 0$ for all $l,k \geqslant 2$, is given by
$$
\alpha_{k,l} := \Delta[ (\Cdot-y)_+^{k+l-2} : b_1,\dots,b_k,c_1,\dots,c_l],\qquad k + l \geqslant 2.
$$
In order to compute $\Delta[(\Cdot -y)^{n-1}_+:a_0,\ldots,a_n] = \alpha_{r,s}$, it suffices therefore to compute the sequence $\alpha_{k,l}$ for
$k+l \geqslant 2$, $k \leqslant r$, $l \leqslant s$, by means of two nested loops, one on $k$, the other on $l$. We detail this computation in
Algorithm~\ref{algo} (see also \cite{Ali73,Var73}).

We can compute $\Delta[(\Cdot -y)^n_-:a_0,\ldots,a_n]$ similarly. Indeed, the same recurrence equation applied to the initial values
$\alpha_{0,l} = 0$ for all $l \geqslant 1$ and $\alpha_{k,0} = 1$ for all $k \geqslant 1$, produces the solution
$$
\alpha_{k,l} := \Delta[ (\Cdot-y)_-^{k+l-1} : b_1,\dots,b_k,c_1,\dots,c_l],\qquad k + l \geqslant 1.
$$

\begin{algorithm}[t!]
\caption{\label{algo} Algorithm for the computation of $\Delta[(\cdot-y)_{+}^{n-1}:a_0,\dots,a_n]$.}
\begin{algorithmic}
\REQUIRE{$n$, $a_0,\dots,a_n$, $y$} \STATE $S \leftarrow 0$, $R \leftarrow 0$ \FOR{$i=0,1,\dots,n$} \IF{$x_i - y \geqslant 0$} \STATE $S
\leftarrow S+1$ \STATE $C_S \leftarrow x_i - y$ \ELSE \STATE $R \leftarrow R+1$ \STATE $B_R \leftarrow x_i - y$ \ENDIF \ENDFOR \STATE $A_0
\leftarrow 0$, $A_1 \leftarrow 1/(C_1 - B_1)$ \COMMENT{Initialization of the unidimensional temporary array of size $S+1$ necessary for the
computation of the divided difference} \FOR{$j=2,\dots,S$} \STATE $A_j \leftarrow - B_1 A_{j-1} / (C_j - B_1)$ \ENDFOR \FOR{$i=2,\dots,R$}
\FOR{$j=1,\dots,S$} \STATE $A_j \leftarrow (C_j A_j - B_i A_{j-1}) / (C_j - B_i)$ \ENDFOR \ENDFOR \STATE {\bf return} $A_R$ \COMMENT{Contains
the value of $\Delta[(\cdot-y)_{+}^{n-1}:a_0,\dots,a_n]$.}
\end{algorithmic}
\end{algorithm}

\begin{exm}
\label{example1} The Choquet integral is frequently used in multicriteria decision aiding, non-additive expected utility theory, or complexity
analysis (see for instance \cite{GraMurSug00} for an overview). For instance, when such an operator is used as an aggregation function in a
given decision making problem, it is very informative for the decision maker to know its distribution. In that context, one of the most natural
{\em a priori}\/ p.d.f.s on $[0,1]^n$ is the standard uniform, which makes the results presented in this section of particular interest.
\begin{figure}[t!]
\begin{center}
\includegraphics*[width=0.7\linewidth]{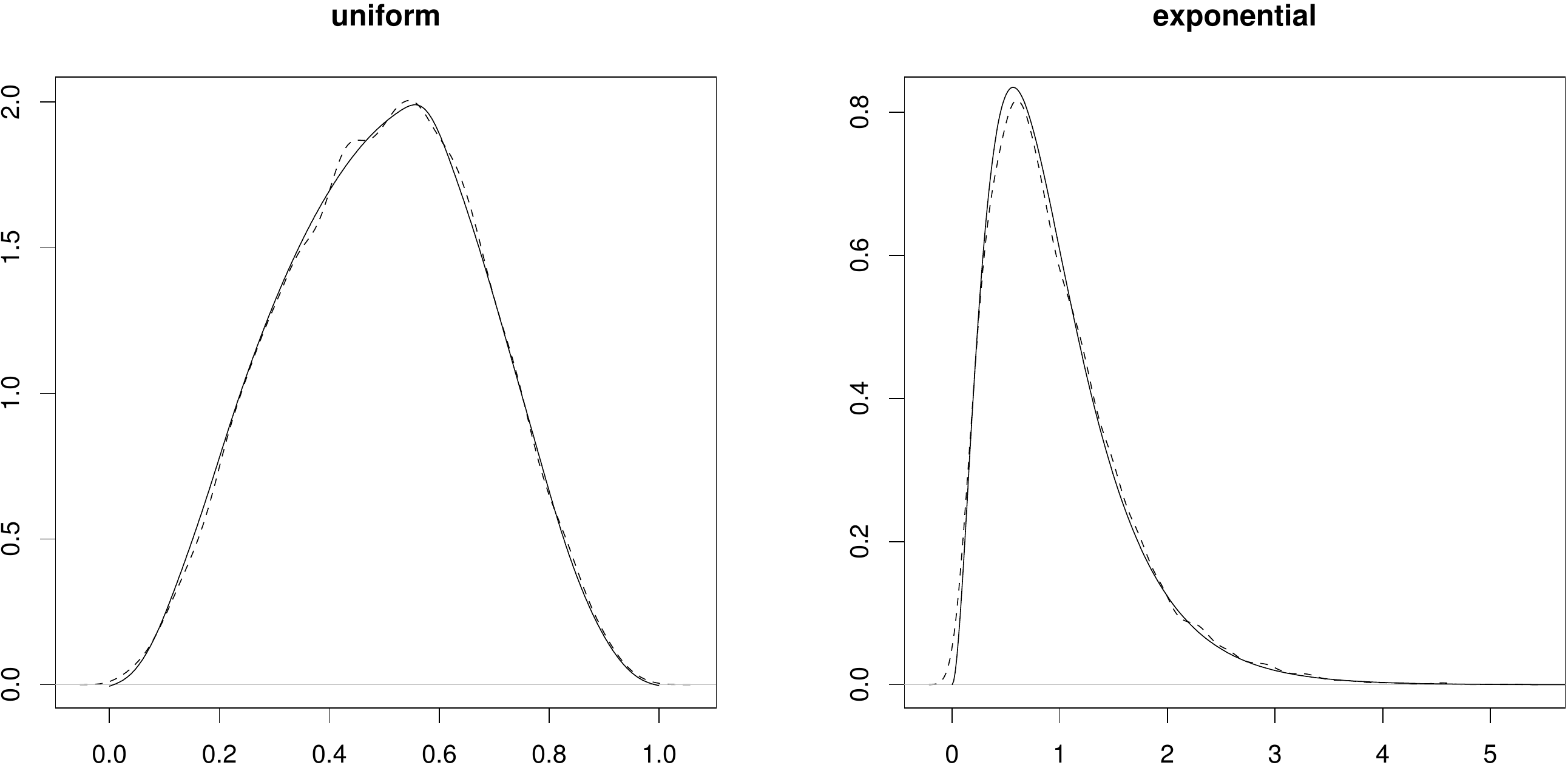}
\caption{\label{fig:cdf_choquet} P.d.f.s of discrete Choquet integral (solid lines) in the standard uniform and standard exponential cases.  The
dotted lines represent the corresponding p.d.f.s estimated by the kernel method from 10~000 randomly generated realizations.}
\end{center}
\end{figure}
Let $\nu$ be the capacity on $N=\{1,2,3\}$ defined by $\nu(\{1\}) = 0.1$, $\nu(\{2\}) = 0.2$, $\nu(\{3\}) = 0.55$, $\nu(\{1,2\}) = 0.7$,
$\nu(\{1,3\}) = 0.8$, $\nu(\{2,3\}) = 0.6$, and $\nu(\{1,2,3\}) = 1$. The p.d.f.\ of the Choquet integral w.r.t.\ $\nu$, which can be computed
through~(\ref{eq:fDens}) and by means of Algorithm~\ref{algo}, is represented in Figure~\ref{fig:cdf_choquet} (left) by the solid line. The
dotted line represents the p.d.f.\ estimated by the kernel method from 10~000 randomly generated realizations  of $U_1,U_2,U_3$ using the R
statistical system~\cite{Rsystem}. The expectation and the standard deviation can also be calculated through~(\ref{ECh_unif})
and~(\ref{ECh2_unif}). We have
$$
\mathbf{E}[Y_\nu]\approx 0.495 \quad\mbox{and}\quad \sqrt{\mathbf{E}[Y_\nu^2]-\mathbf{E}[Y_\nu]^2}\approx 0.183.
$$
The sample mean and the variance of the above mentioned 10~000 realizations of the Choquet integral are
$$
\bar y_\nu \approx 0.497 \quad\mbox{and}\quad s_{y_\nu} \approx 0.183.
$$
\end{exm}

\section{The standard exponential case}

In the standard  exponential case, i.e., when $F(x) = 1 - e^{-x}$, $x\geqslant 0$, the exact distribution of the Choquet integral can be
obtained if the numbers $\{\nu_i^\sigma\}_{i \in N, \sigma \in \frak{S}_n}$ satisfy certain regularity conditions. The result is based on the
following proposition (see  \cite[\S 6.5]{DavNag03} and the references therein).

\begin{prop}\label{prop:ExpCase}
Let $a_1,\dots,a_n \in \R$ and let $X_1,\dots,X_n$ be a random sample drawn from the standard exponential distribution. For any $i \in N$,
define
$$
c_i = \frac{1}{n-i+1} \sum_{j=i}^n a_j.
$$
Then, if $c_i \neq c_k$ whenever $i \neq k$, and $c_i > 0$ for all $i \in N$, the p.d.f.\ of $T=\sum_{i=1}^n a_i X_{i:n}$ is given by
$$
f_T(y) = \sum_{i=1}^n \frac{c_i^{n-2}}{\prod_{k \neq i} (c_i - c_k)} \exp \left( - \frac{y}{c_i} \right).
$$
\end{prop}

The p.d.f.\ $f_\nu(y)$ of the Choquet integral then results from Corollary~\ref{cor:main} and Proposition~\ref{prop:ExpCase}.

\begin{cor}
Assume that, for any $\sigma \in \frak{S}_n$, $\nu_i^\sigma/i \neq \nu_k^\sigma/k$ whenever $i \neq k$, and that $\nu_i^\sigma/i > 0$ for all $i
\in N$. Then
\begin{equation}\label{eq:fDensExp}
f_\nu(y) = \frac{1}{n!} \sum_{\sigma\in\frak{S}_n} \sum_{i=1}^n \frac{(\nu_i^\sigma/i)^{n-2}}{\prod_{k \neq i} (\nu_i^\sigma/i -
\nu_k^\sigma/k)} \exp \left( - \frac{y}{(\nu^\sigma_i/i)} \right).
\end{equation}
\end{cor}

\begin{pf*}{Proof.}
The result is a direct consequence of Corollary~\ref{cor:main}, Proposition~\ref{prop:ExpCase}, and the fact that, for any $\sigma \in
\frak{S}_n$,
$$
\frac{1}{n-i+1} \sum_{j=i}^n p_{n-j+1}^{\nu,\sigma} = \frac{\nu_{n-i+1}^\sigma}{n-i+1}, \qquad i \in N.\qed
$$
\qed \end{pf*}

The first two moments of the order statistics in the standard exponential case are given (see e.g.\ \cite[p.~52]{DavNag03}) by
\begin{equation}
\label{EXi_exp} \mathbf{E}[X_{i:n}] = \sum_{k=n-i+1}^n \frac{1}{k},
\end{equation}
and, if $i < j$,
\begin{equation}
\label{EXiXj_exp} \mathbf{E}[X_{i:n} X_{j:n}]  - \mathbf{E}[X_{i:n}] \mathbf{E}[X_{j:n}] = \mathbf{E}[X_{i:n}^2]  -
\mathbf{E}[X_{i:n}]^2=\sum_{k=n-i+1}^n \frac{1}{k^2}.
\end{equation}
Used in combination with~(\ref{eq:expect}) and~(\ref{eq:expect_square}), these expressions enable us to obtain the first two raw moments of the
Choquet integral.

\begin{exm}
\label{example2} Consider again the capacity given in Example~\ref{example1} and assume now that $X_1,X_2,X_3$ is a random sample from the
standard exponential distribution. The p.d.f.\ of the Choquet integral w.r.t.\ $\nu$, which can be computed by means of~(\ref{eq:fDensExp}), is
represented in Figure~\ref{fig:cdf_choquet} (right) by the solid line. The dotted line represents the p.d.f.\ estimated by the kernel method
from 10~000 randomly generated realizations.

Combining~(\ref{EXi_exp}) and~(\ref{EXiXj_exp}) with~(\ref{eq:expect}) and~(\ref{eq:expect_square}), we obtain the following values:
$$
\mathbf{E}[Y_\nu]\approx 0.963 \quad\mbox{and}\quad \sqrt{\mathbf{E}[Y_\nu^2]-\mathbf{E}[Y_\nu]^2}\approx 0.624.
$$
The sample mean and the variance of the above mentioned 10~000 realizations of the Choquet integral are
$$
\bar y_\nu \approx 0.964 \quad\mbox{and}\quad s_{y_\nu} \approx 0.630.
$$

\end{exm}
\section{Approximations of the moments}
\label{approx_moments}

When $F$ is neither the standard uniform, nor the standard exponential c.d.f., but $F^{-1}$ and its derivatives can be easily computed, one can
obtain approximations of the moments of order statistics, and therefore of those of the Choquet integral, using the approach initially proposed
by David and Johnson~\cite{DavJoh54}.

Let $U_1,\dots,U_n$ be a random sample from the standard uniform distribution. The product moments of the corresponding order statistics are then given by the following formula: 
\begin{equation}
\label{moments_uniform} \mathbf{E} \Big[ \prod_{j=1}^l U_{i_j:n}^{m_j} \Big] = \frac{n!}{\left(n + \sum_{j=1}^l m_j \right)! } \prod_{j=1}^l
\frac{(i_j + m_1 + \dots + m_j - 1)!}{(i_j + m_1 + \dots + m_{j-1} - 1)!},
\end{equation}
where $1 \leqslant i_1<\cdots< i_l \leqslant n$. Now, it is well known that the c.d.f.\ of $X_{i:n}$ is given by
$$
\mathrm{Pr}[X_{i:n} \leqslant x] = \sum_{j=i}^n \binom{n}{j} F^j(x) [1 - F(x)]^{n-j}.
$$
It immediately follows that
$$
\mathrm{Pr}[F^{-1}(U_{i:n}) \leqslant x] = \mathrm{Pr}[U_{i:n} \leqslant F(x)] = \mathrm{Pr}[X_{i:n} \leqslant x],
$$
i.e., that $F^{-1}(U_{i:n})$ and $X_{i:n}$ are equal in distribution.

Starting from this distributional equality, David and Johnson~\cite{DavJoh54} expanded $F^{-1}(U_{i:n})$ in a Taylor series around the point
$\mathbf{E}[U_{i:n}] = i/(n+1)$ in order to obtain approximations of product moments of non-uniform order statistics. Setting $r_i := i/(n+1)$,
$G := F^{-1}$, $G_i := G(r_i)$, $G^{(1)}_i := G^{(1)}(r_i)$, etc., we have
$$
X_{i:n} = G_i + (U_{i:n} - r_i)G^{(1)}_i + \frac{1}{2} (U_{i:n} - r_i)^2 G^{(2)}_i + \frac{1}{6} (U_{i:n} - r_i)^3 G^{(3)}_i + \dots
$$
Setting $s_i := 1 - r_i$, taking the expectation of the previous expression and using~(\ref{moments_uniform}), the following approximation for
the expectation of $X_{i:n}$ can be obtained to order $(n+2)^{-2}$:
\begin{equation}
\label{EXi} \mathbf{E}[X_{i:n}] \approx G_i + \frac{r_i s_i}{2(n+2)} G^{(2)}_i + \frac{r_i s_i}{(n+2)^2} \left[ \frac{1}{3} (s_i - r_i)
G^{(3)}_i + \frac{1}{8} r_i s_i G^{(4)}_i \right].
\end{equation}
%
Similarly, for the first product moment, we have
\begin{align}
\nonumber \mathbf{E}[X_{i:n}X_{j:n}] \approx& \, G_i G_j + \frac{r_i s_j}{n+2} G^{(1)}_i G^{(1)}_j + \frac{r_i s_i}{2 (n+2)}G_j G^{(2)}_i + \frac{r_j s_j}{2 (n+2)} G_i G^{(2)}_j  \\
\nonumber &+ \frac{r_i s_j}{(n+2)^2} \Big[ (s_i - r_i) G^{(2)}_i G^{(1)}_j + (s_j - r_j) G^{(1)}_i G^{(2)}_j + \frac{1}{2} r_i s_i G^{(3)}_i G^{(1)}_j  \\
\nonumber &+ \frac{1}{2} r_j s_j G^{(1)}_i G^{(3)}_j  +  \frac{1}{2} r_i s_j G^{(2)}_i G^{(2)}_j \Big] + \frac{r_i r_j s_i s_j}{4
(n+2)^2}G^{(2)}_i G^{(2)}_j
\\ \nonumber &+ \frac{r_i s_i G_j}{(n+2)^2} \Big[ \frac{1}{8} r_i s_i G^{(4)}_i + \frac{1}{3} (s_i -r_i) G^{(2)}_i \Big] \\
\label{EXiXj} &+ \frac{r_j s_j G_i}{(n+2)^2} \Big[\frac{1}{8} r_j s_j G^{(4)}_j + \frac{1}{3} (s_j - r_j) G^{(2)}_j \Big].
\end{align}
The accuracy of the above approximations is discussed in \cite[\S 4.6]{DavNag03}. Note that Childs and Balakrishnan~\cite{ChiBal02} have
recently proposed MAPLE routines facilitating the computations and permitting the inclusion of higher order terms.

As already mentioned, the previous expressions are useful only if $G:=F^{-1}$ and its derivatives can be easily computed. This is the case for
instance when $F$ is the standard normal c.d.f. Indeed, there exist algorithms that enable an accurate computation of $F^{-1}$ and it can be
verified (see e.g.\ \cite[p 85]{DavNag03}) that $G^{(1)}=(f \circ G)^{-1}$,
$$
G^{(2)} = \frac{G}{f^2 \circ G}, \qquad G^{(3)} = \frac{1+2G^2}{f^3 \circ G} \qquad \mbox{and} \qquad G^{(4)}=\frac{G(7+6G^2)}{f^4 \circ G},
$$
where $f := F^{(1)}$.

From a practical perspective, in order to obtain a better accuracy for $\mathbf{E}[X_{i:n}]$ and $\mathbf{E}[X_{i:n}X_{j:n}]$ in the standard
normal case, one can use the expressions obtained to order $(n+2)^{-3}$ in \cite{DavJoh54} and recalled in \cite{ChiBal02}. We do not reproduce
these expressions here as they are very long. We provide however the expressions of $G^{(5)}$ and $G^{(6)}$ required for computing them:
$$
G^{(5)} = \frac{7+G^2(46+24 G^2)}{f^5 \circ G} \qquad \mbox{and} \qquad G^{(6)}=\frac{G(127+326G^2 + 96G^4)}{f^6 \circ G}.
$$

\begin{exm}
\label{example3} Consider again the capacity given in Example~\ref{example1} and assume now that the decision maker wants the standard normal as
{\em a priori} p.d.f. Combining~(\ref{EXi}) and~(\ref{EXiXj}) with~(\ref{eq:expect}) and~(\ref{eq:expect_square}), we obtain the following
approximate values:
$$
\mathbf{E}[Y_\nu]\approx -0.014 \quad\mbox{and}\quad \sqrt{\mathbf{E}[Y_\nu^2]-\mathbf{E}[Y_\nu]^2}\approx 0.615.
$$
For comparison, the sample mean and the variance of 10~000 independent realizations of the corresponding Choquet integral are
$$
\bar y_\nu \approx -0.013 \quad\mbox{and}\quad s_{y_\nu} \approx 0.620.
$$

\end{exm}


\section{Asymptotic distribution of the Choquet integral}

Conditions under which a linear combination of order statistics is asymptotically normal have been extensively studied in the statistical
literature. A good synthesis on the subject is given in \cite[\S 11.4]{DavNag03}. Provided some regularity conditions are satisfied, typically
on $\nu$ and $F$ in the context under consideration, the existing theoretical results, combined with Proposition~\ref{prop:main}, practically
imply that, for large $n$, $Y_\nu$ is approximately distributed as a mixture of $n!$ normals
$N(\mathbf{E}[Y_\nu^\sigma],\mathbf{V}[Y_\nu^\sigma])$, $\sigma \in \frak{S}_n$, each weighted by $\frac{1}{n!}$.

From a practical perspective, the most useful result seems to be that of Stigler~\cite{Sti74}. For any $\sigma \in \frak{S}_n$, let
$J^{\nu,\sigma}$ be a real function on $[0,1]$ such that $J^{\nu,\sigma}(i/n) = n p^{\nu,\sigma}_{n-i+1}$. Then, $Y_\nu^\sigma$ can be rewritten
as
$$
Y_{\nu,n}^\sigma = \frac{1}{n} \sum_{i=1}^n J^{\nu,\sigma} \left( \frac{i}{n} \right) X_{i:n},
$$
where the subscript $n$ in $Y_{\nu,n}^\sigma$ is added to emphasize dependence on the sample. Furthermore, let
$$
\alpha(J^{\nu,\sigma},F) := \int_{-\infty}^\infty x J^{\nu,\sigma}[F(x)] dF(x)
 = \int_0^1 J^{\nu,\sigma}(u) F^{-1}(u) \mathrm{d}u,
$$
and
\begin{align*}
\beta^2(J^{\nu,\sigma},F):&=2 \int_{-\infty<x<y<+\infty} J^{\nu,\sigma}(F(x))J^{\nu,\sigma}(F(y)) F(x)(1-F(y)) \mathrm{d} x \mathrm{d} y \\
&=2 \int_{0<u<v<1} J^{\nu,\sigma}(u)J^{\nu,\sigma}(v) u(1-v) \mathrm{d}F^{-1}(u) \mathrm{d}F^{-1}(v).
\end{align*}
Then, Stigler's results~\cite[Theorems~2 and~3]{Sti74} (see also \cite[Theorem~11.4]{DavNag03}) state that, if $F$ has a finite variance and if
$J^{\nu,\sigma}$ is bounded and continuous almost everywhere w.r.t.\ $F^{-1}$, one has
$$
\lim_{n \rightarrow \infty} \mathbf{E}[Y_{\nu,n}^\sigma] = \alpha(J^{\nu,\sigma},F), \qquad \lim_{n \rightarrow \infty} n
\mathbf{V}[Y_{\nu,n}^\sigma] = \beta^2(J^{\nu,\sigma},F),
$$
and, if additionally $\beta^2(J^{\nu,\sigma},F) > 0$,
$$
\frac{Y_{\nu,n}^\sigma - \mathbf{E}[Y_{\nu,n}^\sigma]}{\sqrt{\mathbf{V}[Y_{\nu,n}^\sigma]}} \rightarrow_{d} N(0,1) \quad \mbox{as} \quad n
\rightarrow \infty.
$$

\begin{exm}
To illustrate the applicability of these results, consider the following game $\nu$ on $N$ defined by
$$
\nu(S) = \sum_{j=1}^{|S|} \frac{1}{n} \left( \frac{n-j+1}{n} \right)^a, \qquad \forall S \subseteq N,
$$
where $a$ is a strictly positive real number. We then have
$$
p_i^{\nu,\sigma} = \frac{1}{n} \left( \frac{n-i+1}{n} \right)^a, \qquad \forall i \in N, \qquad \forall \sigma \in \frak{S}_n.
$$
As the coefficients $p_i^{\nu,\sigma}$ do not depend on $\sigma$, the corresponding Choquet integral is merely a linear combination of order
statistics. Note however that the game $\nu$ is by no means additive. Next, define $J^{\nu,\sigma}(x) := x^a$, for all $x \in [0,1]$. Then,
clearly, $J^{\nu,\sigma}(i/n) = n p_{n-i+1}^{\nu,\sigma}$ for all $i \in N$.

In order to simplify the calculations, assume furthermore that $F$ is the standard uniform c.d.f.\ and that $a=2$. Then, $J^{\nu,\sigma}$ is
clearly bounded and continuous almost everywhere w.r.t.\ $F^{-1}$ and we have $\alpha(J^{\nu,\sigma},F) = 1/4$ and $\beta^2(J^{\nu,\sigma},F) =
1/112$.

\begin{figure}[t!]
\begin{center}
\includegraphics*[width=0.9\linewidth]{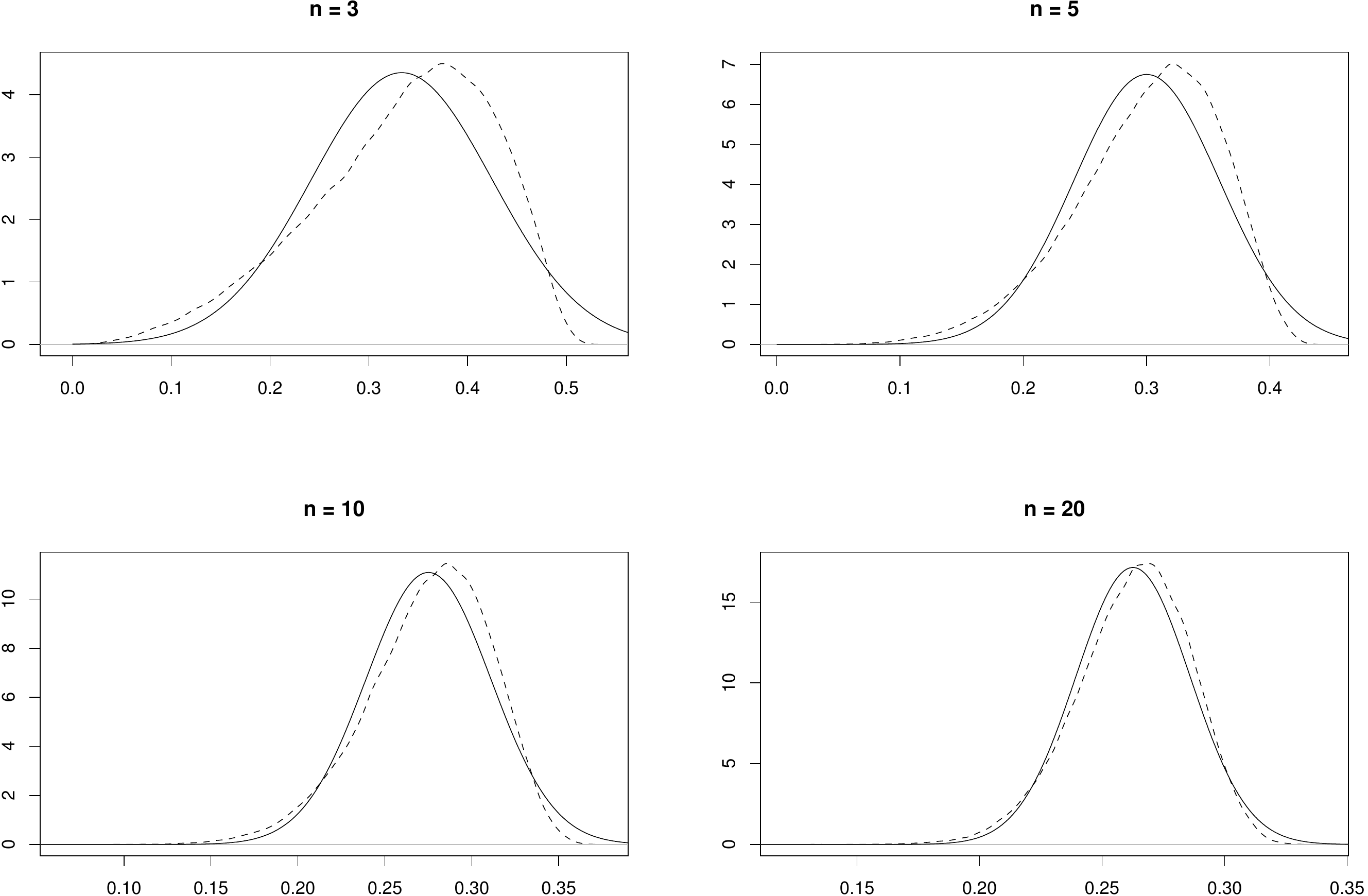}
\caption{\label{fig:asym_choquet} Approximations of the p.d.f.s of discrete Choquet integral by mixtures of normals (solid lines) for $n=3,5,10$
and 20. The dotted lines represent the corresponding p.d.f.s estimated by the kernel method from 10~000 randomly generated realizations.}
\end{center}
\end{figure}

The dotted lines in Figure~\ref{fig:asym_choquet} represent the p.d.f.\ of the Choquet integral w.r.t.\ $\nu$ estimated by the kernel method
from 10~000 randomly generated realizations for $n=3,5,10$ and 20. The solid lines represent the normal p.d.f.s
$N(\mathbf{E}[Y_{\nu,n}^\sigma],\mathbf{V}[Y_{\nu,n}^\sigma])$, where $\mathbf{E}[Y_{\nu,n}^\sigma]$ and $\mathbf{V}[Y_{\nu,n}^\sigma]$ are
computed by means of~(\ref{ECh_unif}) and~(\ref{ECh2_unif}).
\end{exm}

From the previous example, it clearly appears that one strong prerequisite before being able to apply the previous theoretical results is the
knowledge of the expression of the game $\nu$ in terms of $n$. In practical applications of aggregation operators, this is rarely the case as
$\nu$ is usually determined for some fixed $n$ from learning data (see e.g.\ \cite{GraKojMey08}). It follows that in such situations the above
theoretical conditions cannot be rigorously verified.

In informal terms, Stigler~\cite{Sti74} states that a linear combination of order statistics is likely to be asymptotically normally distributed
if the extremal order statistics do not contribute ``too much'', which is satisfied is the weights are ``smooth'' and ``bounded''. When dealing
with a Choquet integral, several numerical indices could be computed to assess whether the operator behaves in a too {\em conjunctive}
(minimum-like) or too {\em disjunctive} (maximum-like) way. One such index is the degree of {\em orness} studied in \cite{Mar00b,Mar04}.

\begin{figure}[t!]
\begin{center}
\includegraphics*[width=1\linewidth]{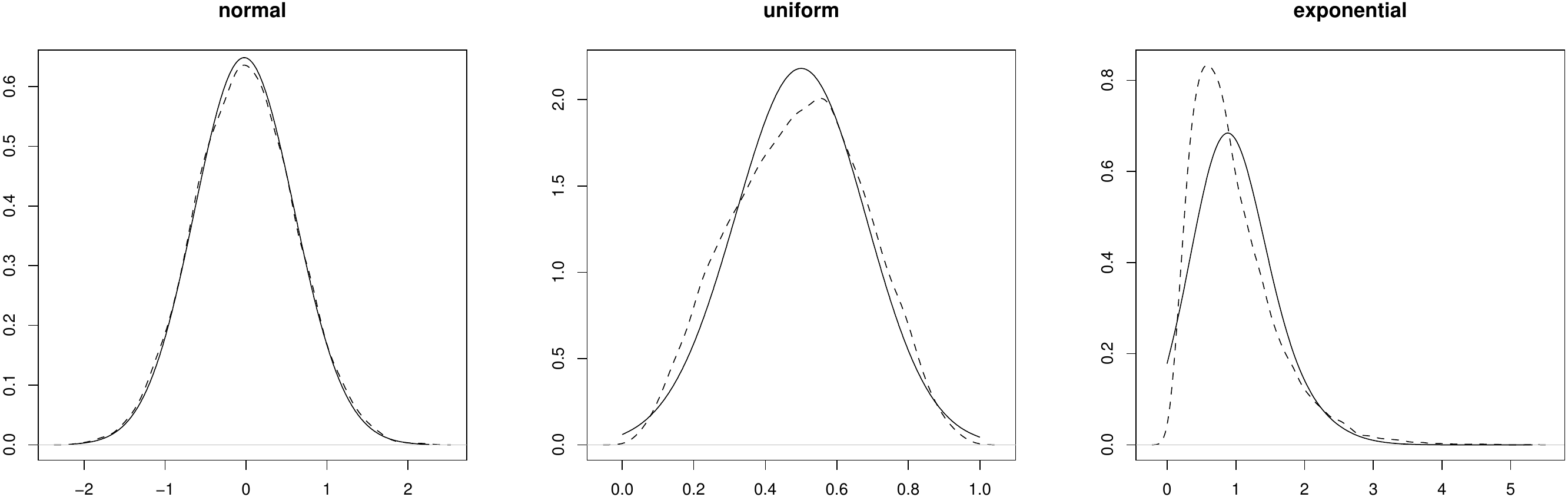}
\caption{\label{fig:asym_choquet_2} Approximations of the p.d.f.s of discrete Choquet integrals by mixtures of normals (solid lines) in the
standard normal, standard uniform and standard exponential cases. The dotted lines represent the corresponding p.d.f.s estimated by the kernel
method from 10~000 randomly generated realizations.}
\end{center}
\end{figure}

\begin{exm}
\label{example4} Consider again the capacity given in Example~\ref{example1}. The degree of orness of this capacity, computed using the {\tt
kappalab} R package, is 0.49, which indicates a fairly neutral (slightly conjunctive) behavior.  The solid lines in
Figure~\ref{fig:asym_choquet_2} represent the mixtures of $3!=6$ normals in the standard normal, standard uniform and standard exponential cases
as possible approximations of the p.d.f.\ of the corresponding Choquet integral. As previously, the dotted lines represent the  p.d.f.s
estimated by the kernel method from 10~000 randomly generated realizations. As one can see, the approximation is very good in the standard
normal case, may be considered as acceptable in the standard uniform case, and poor in the exponential case. Provided considering such a
approximation is valid (which, as discussed above, cannot be verified), one could argue that the poor results in the exponential case are due to
the too low value of $n(=3)$. Although such low values for $n$ make no sense in statistics, in multicriteria decision aiding for instance, they
are quite common. In fact, in practical decision problems involving aggregation operators, the value of $n$ is very rarely greater than 10.
\end{exm}

\end{document}